\documentclass[11pt,letterpaper]{article}

\usepackage{hyperref}

\usepackage{amssymb,amsmath}

\usepackage{epsf}
\usepackage{graphicx,amssymb,amsmath,color,latexsym} 

\usepackage{xcolor,wasysym,floatflt}

\newtheorem{theorem}{Theorem}

\newtheorem{rem}{Remark}
\newtheorem{pro}{Proposition}

\newcommand{\eps}{\varepsilon}

\usepackage[none]{hyphenat} 

\makeatletter
\newcommand\footnoteref[1]{\protected@xdef\@thefnmark{\ref{#1}}\@footnotemark}
\makeatother
\def\qed{\hfill $\square$}

\begin{document}

\begin{center}
   \Large \textbf{Bifurcations of finite-time stable limit cycles from focus boundary equilibria in impacting systems, Fillippov systems and sweeping processes}  
 \vskip 0.3cm \large Oleg Makarenkov, Lakmi Niwanthi Wadippuli
  \vskip.1cm \small Department of Mathematical Sciences, The University of Texas at Dallas, 800 West Campbell Road 
  Richardson, Texas 75080, USA
\end{center}

\noindent Abstract: We establish a theorem on bifurcation of limit cycles from a focus boundary equilibrium of an impacting system, which is universally applicable to prove bifurcation of limit cycles from focus boundary equilibria in other types of piecewise-smooth systems, such as Filippov systems  and sweeping processes. Specifically, we assume that one of the subsystems of the piecewise-smooth system under consideration admits a focus equilibrium that lie on the switching manifold at the bifurcation value of the parameter. In each of the three cases, we derive a linearized system which is capable to conclude about the occurrence of a finite-time stable limit cycle from the above-mentioned focus equilibrium when the parameter crosses the bifurcation value. Examples illustrate how conditions of our theorems lead to closed-form formulas for the coefficients of the linearized system.

%\begin{figure}[h]\center \noindent\includegraphics[scale=0.82]{P.pdf}\end{figure}

\section{Introduction}

\noindent Unfolding of a singular equilibrium of a vector field on a boundary of a smooth manifold is a classical problem of the theory of differential equations that goes back to Vishik  \cite{vishik} and Arnold \cite{arnold}. 

%\noindent Classical Hopf bifurcation remains one of the central mechanisms used by applied scientists to establish the occurrence of smooth limit cycles in nonlinear mathematical models. However, limit cycles of many important dynamical processes are fundamentally nonsmooth. Some outstanding examples are robotic locomotion \cite{grizzle}, dc-dc power converting \cite{schild,ian}, anti-lock braking \cite{astolfi,oleg}, neuron firing \cite{izh,coombes}. That is why the theory of bifurcation of nonsmooth limit cycles has been intensively developing lately. As in the classical Hopf bifurcation, nonsmooth limit cycles occur from an equilibrium that gets singular in a certain nonsmooth sense. Specifically, the most typical scenarios that leads to the occurrence of a nonsmooth limit cycle is when an equilibrium of a piecewise smooth dynamical system approaches the switching manifold under varying parameters and attempts to cross it. Such a scenarios is known  as border-collision bifurcation in nonsmooth dynamical systems literature.

\vskip0.2cm

\noindent In the case where the boundary of a smooth manifold is a switching manifold separating two smooth differential equations, the main breakthrough is due to Filippov \cite{filippov}, who offered a formula to define the flow of the full (i.e. piecewise smooth) system of differential equations on the switching manifold (called {\it sliding flow}). In particular, Filippov observed \cite[\S~19]{filippov}  that a focus equilibrium of a smooth subsystem of a piecewise smooth planar system of differential equations may produce a limit cycle after such an equilibrium collides with the switching manifold under varying parameters. In this way Filippov paved a route to such an  analogue of the classical Hopf bifurcation that is capable to provide limit cycles that lack smoothness (with multiple applications to e.g. mechanical systems with dry friction \cite{oleg}).

\vskip0.2cm

\noindent The problem of bifurcation of limit cycles from focus boundary equilibria of Filippov systems has been intensively refined lately, see e.g. Kuznetsov et al, Guardia et al, Hogan et al, Glendinning. Specifically, if a  Filippov system 
\begin{equation}\label{0}
\left( \begin{array}{c}\dot x\\ \dot y\end{array}\right)=\left( \begin{array}{c}f^i(x,y,\eps)\\ g^i(x,y,\eps)\end{array}\right),\quad i=\left\{\begin{array}{l}
i=+1,\ {\rm if}\ H(x,y)>0,\\
i=-1,\ {\rm if}\ H(x,y)<0,
\end{array}\right.
\end{equation}
where $f^i,$ $g^i$ and $H$ are smooth functions,
admits a focus equilibrium $(x_\eps,y_\eps)\to (x_0,y_0)$ as $\eps\to 0$ with $H(x_0,y_0)=0$, then the available theory (as it appears e.g. in \cite{gle}) provides a change of the variables that brings (\ref{0}) near $(x_0,y_0)$  
to the normal form 
\begin{eqnarray}
   &&\left(\begin{array}{c} \dot x \\
   \dot y\end{array}\right)=\left(\begin{array}{cc} a & -b \\ b & a \end{array}\right)\left(\begin{array}{c} x \\ y-\eps\end{array}\right)+\mbox{smaller nonlinear terms},\quad \mbox{if} \ y>0,\label{lin1_}\\
  &&\left(\begin{array}{c} \dot x \\
   \dot y\end{array}\right)=\left(\begin{array}{c} m \\ 1\end{array}\right)+\mbox{smaller nonlinear terms},\quad \mbox{if} \ y<0.\label{lin2_}
\end{eqnarray}
One of the conclusions of Glendenning \cite{gle} and  Kuznetsov et al \cite{kuz} relate the property of the form 
\begin{equation}\label{parameters}
\mbox{either }m>\dfrac{a}{b}\quad \mbox{ or }\quad m<\dfrac{a}{b}\mbox{ and }\dfrac{1}{\eps}P\left(\dfrac{a}{b}\eps,\eps\right)<\dfrac{am+b}{bm-a}, 
\end{equation} 
to the existence of cycles in the linear part of system (\ref{lin1_})-(\ref{lin2_}) for $\eps>0.$ Here $x\to P(x,\eps)$ is the Poincar\'e map of (\ref{lin1_}) induced by the cross-section $y=0$. And the purpose of the second inequality of (\ref{parameters}) is to avoid the presence of stationary points of the sliding flow between $\dfrac{a}{b}\eps$ and $\dfrac{1}{\eps}P\left(\dfrac{a}{b}\eps,\eps\right).$

\vskip0.2cm

\noindent Much less is known in the case where a focus equilibrium collides with the boundary of a completely inelastic unilateral constraint, which formulates as a differential inclusion (see e.g. \cite{kunze,edm})
$$
   \left(\begin{array}{c}\dot x\\ \dot y\end{array}\right)\in -N_{C(\eps)}(x,y)+\left(\begin{array}{c}  f(x,y) \\ g(x,y)\end{array}\right),
$$
where $N_C(x,y)$ is Clarke's normal cone to $C$ at $(x,y)$,
and is known as sweeping process. Sweeping processes is a standard tool to describe the evolution of elastoplastic systems \cite{bas}, that currently gain attention also in the context of crown motion modeling \cite{mord}. Here the very concept of a stationary point of a sliding flow has been defined just recently (in slightly different terminology) and no any results about bifurcations from boundary equilibria are currently available, that was the main motivation of our work.

\vskip0.2cm

\noindent In this paper we offer a unified theorem on bifurcation of limit cycles from a boundary equilibrium  of a  hybrid system
\begin{eqnarray}
&\left( \begin{array}{c}\dot x\\ \dot y\end{array}\right)=\left( \begin{array}{c}f(x,y)\\ g(x,y)\end{array}\right),& H(x,y,\eps)<0,\label{1a} \\
& (x,y)\mapsto \left(A(\eps),B(\eps)\right), &  H(x,y,\eps)=0,\label{1b}
\end{eqnarray} which is capable to predict the occurrence of limit cycles  in Filippov systems and sweeping processes alike. Compared to the above-mentioned results about bifurcation of limit cycles in Filippov systems, our result implies the occurrence of a cycle in the initial nonlinear system (\ref{0}), rather than in its linearization given by (\ref{lin1_})-(\ref{lin2_}). The linear system of (\ref{lin1_})-(\ref{lin2_}) is considered as an example in which case we get same condition (\ref{parameters}). Note, following di Bernardo et al \cite{olivar}, a different equivalent strategy can be taken where bifurcation results in both impact systems and sweeping processes are derived from a general result for Filippov systems. The later strategy has been also offered earlier by Zhuravlev \cite{zh} and Ivanov \cite{iv}.

\vskip0.2cm

\noindent The paper is organized as follows. In Section~2 we prove our main result (Theorem~\ref{thm1}) about bifurcation of limit cycles in hybrid system (\ref{1a})-(\ref{1b}) from the origin 
which is a focus equilibrium of subsystem (\ref{1a}). We consider parameter-independent vector fields in (\ref{1a}), but  rather assume that the switching manifold is a function of the parameter $\eps$ and that the origin belongs to the switching boundary at $\eps=0$ (i.e. that $H(0,0,0)=0$). To illustrate Theorem~\ref{thm1} a simple resonate-and-fire neuron model from Izhikevich \cite{izh1} is considered. In di Bernardo et al \cite{olivar}, the analysis of bifurcations of limit cycles from a boundary focus equilibrium in impact system (\ref{1a})-(\ref{1b}) is converted into the analysis of the respective bifurcations in Filippov systems, but  the approach of \cite{olivar} uses state-dependence of the impact law (\ref{1b}) in an essential way.

\vskip0.2cm

\noindent Section~3 shows (Theorem~\ref{thm2}) that bifurcation of limit cycles in Filippov system of type (\ref{0}) from a boundary focus equilibrium, can be obtained as a corollary of Theorem~\ref{thm1}. We note that throughout Section~3 we assume that vector fields in (\ref{0}) don't depend on $\eps,$ but $H$ does, which is equivalent to the setting (\ref{0}). The linear part of (\ref{lin1_})-(\ref{lin2_}) is considered in Section~3 as a benchmark to illustrate Theorem~\ref{thm2}. Here we also enhance the known formula
 (\ref{parameters}) by  deriving a closed-form expression for the last inequality of (\ref{parameters}), that allows us to plot (\ref{parameters}) in the $\left(\dfrac{a}{b},m\right)$-coordinate plane (Fig.~\ref{fig1}).
 
 \vskip0.2cm
 
 \noindent An application of Theorem~\ref{thm1} to sweeping processes is given in Section~4. The properties similar to those of the  Filippov sliding vector field are established for sliding along the boundary $\partial C(\eps)$ of the unilateral constraint $C(\eps)$ in 
 Proposition~\ref{pro2} of Section~4. In particular, formula (\ref{ssp}) introduces an equation of sliding along $\partial C(\eps)$ and formula (\ref{construction}) gives an equation for stationary point of sliding motion. Based on the properties discovered in Proposition~\ref{pro2}, Theorem~\ref{thm2} establishes bifurcation of a finite-time stable limit cycle as $\partial C(\eps)$ collides with a focus equilibrium of the vector field $\left(\begin{array}{c} f(x,y)\\ g(x,t)\end{array}\right)$ of perturbed sweeping process.

\section{Impacting systems}
The change of the variables
$$
  \left(\begin{array}{c}
  u(t)\\ v(t)\end{array}\right)=\dfrac{1}{\eps}\left( \begin{array}{c} x\\  y\end{array}\right)
$$
brings (\ref{1a})-(\ref{1b}) to the form
\begin{eqnarray}
&\left( \begin{array}{c}\dot u\\ \dot v\end{array}\right)=\dfrac{1}{\eps}\left( \begin{array}{c}f(\eps u,\eps v)\\ g(\eps u,\eps v)\end{array}\right),& H(\eps u,\eps v,\eps)<0,\label{2a} \\
& (u,v)\mapsto \dfrac{1}{\eps}\left(A(\eps),B(\eps)\right), &  H(\eps u,\eps v,\eps)=0.\label{2b}
\end{eqnarray}
We identify $(u,v)$ and $(u,v)^T$ when it doesn't lead to a confusion.
Along with system (\ref{2a})-(\ref{2b}) we consider the following reduced system
\begin{eqnarray}
& \left( \begin{array}{c}\dot u\\ \dot v\end{array}\right)=\left(\begin{array}{cc}f_x(0)  & f_y(0) \\g_x(0) & g_y(0) \end{array}\right) \left( \begin{array}{c} u\\ v\end{array}\right), & {\rm if}\ H_y(0)v+H_\eps(0)<0\label{3a}\\
&u\mapsto A'(0), & {\rm if}\ v=-\dfrac{H_\eps(0)}{H_y(0)}.\label{3b}
\end{eqnarray}

\begin{theorem}\label{thm1} Assume  that the equilibrium of (\ref{1a}) collides with the switching manifolds when $\eps=0$, i.e. $f(0)=g(0)=H(0)=0$. Assume that the coordinates are rotated in such a way that 
 $H_x(0)=0$ and $H_y(0)\not=0$. Assume that the vector field of (\ref{1a}) is tangent to the switching manifold at $(A(\eps),B(\eps))$, i.e., for all $\eps>0,$
\begin{equation}\label{as1rev}
  H(A(\eps),B(\eps),\eps)=0
  %\ \ \mbox{and}\ \  \left(H_x(A(\eps),B(\eps),\eps),H_y(A(\eps),B(\eps),\eps)\right)\left(\begin{array}{c}f(A(\eps),B(\eps)) \\ g(A(\eps),B(\eps))\end{array}\right)=0.
\end{equation}
Assume that the reduced system (\ref{3a})-(\ref{3b}) admits a cycle $(u_0(t),v_0(t))$  with the initial condition $(u_0(0),v_0(0))=\left(A'(0),-\dfrac{H_\eps(0)}{H_y(0)}\right)$
of exactly one impact per period.
  Let $T_0$ be the period of the cycle. If
  \begin{equation}\label{as2rev}
     %u_0(T_0)\not=A'(0),\quad 
   g_x(0)u_0(T_0)-g_y(0)\dfrac{H_\eps(0)}{H_y(0)}\not=0,\quad
     H_\eps(0)\not=0,
  \end{equation}
  then, for all $\eps>0$ sufficiently small, the impacting system (\ref{1a})-(\ref{1b}) admits a finite-time stable limit cycle $(x_\eps(t),y_\eps(t))$  with the initial condition  $(x_\eps(0),y_\eps(0))=(A(\eps),B(\eps)).$ Specifically, there exists $T_\eps\to T_0$ as $\eps\to 0$ such that $H(x_\eps(T_\eps),y_\eps(T_\eps))=0,$ for all $\eps>0$ sufficiently small. 
  \end{theorem}

\noindent {\bf Proof.} Let %$t\mapsto \left(\begin{array}{c}X(t,x,y)\\  Y(t,x,y)\end{array}\right)$ and 
$t\mapsto \left(\begin{array}{c}
U(t,u,v,\eps)\\  V(t,u,v,\eps)\end{array}\right)$ be the general solution of  system 
%of (\ref{1a}) and 
(\ref{2a}). 
 Introduce
$$
   F(T,\eps)=\dfrac{1}{\eps}H\left(\eps U\left(T,\dfrac{A(\eps)}{\eps},\dfrac{B(\eps)}{\eps},\eps\right),\eps V\left(T,\dfrac{A(\eps)}{\eps},\dfrac{B(\eps)}{\eps},\eps\right),\eps\right).
$$
Computing $F(T,0)$ we get
$$
   F(T,0)=H_y(0)V(T,A'(0),B'(0))+H_\eps(0).
$$
The value of $B'(0)$ can be found from (\ref{as1rev}) as 
$$
  B'(0)=-\dfrac{H_\eps(0)}{H_y(0)}.
$$
Therefore,  $F(T_0,0)=0,$ and   since
$$
   F_t(0,T_0)=H_y(0)V_t(T_0,A'(0),B'(0),0)=H_y(0)(g_x(0),g_y(0))\left(\begin{array}{c}
   U(T_0,A'(0),B'(0),0)\\ V(T_0,A'(0),B'(0),0)
   \end{array}\right),
$$ we have 
$F_t(T_0,0)\not=0$ by the first assumption of (\ref{as2rev}).
 Therefore, the existence of $T_\eps$ such that $F(T_\eps,\eps)=0$ follows by applying the Implicit Function Theorem, which in turn implies that $(x_\eps(t),y_\eps(t))$ is a cycle of (\ref{1a})-(\ref{1b}).

\vskip0.2cm

\noindent To establish finite-time stability of $(x_\eps(t),y_\eps(t))$ we have to prove that $(x_\eps(t),y_\eps(t))$ reaches the switching manifold $L=\{(x,y)\in\mathbb{R}^2:H(x,y,\eps)$ transversally (see also \cite[Proposition~1]{oleg}). In other words, we have to show that 
$$
  \phi(\eps)=\left(H_x(u(T_\eps),v(T_\eps),\eps),H_y(u(T_\eps),v(T_\eps)\right)\left(\begin{array}{c} u(T_\eps) \\ v(T_\eps)\end{array}\right)
$$
doesn't vanish for all $\eps>0$. Indeed, we have $\phi(0)=B'(0)H_y(0)=-H_\eps(0)\not=0$ by the second assumption of (\ref{as2rev}).

\vskip0.2cm

\noindent The proof of the theorem is complete.\qed

\vskip0.2cm

\noindent As an example we consider the following nonlinear model of a resonate-and-fire neuron from  Izhikevich \cite{izh1}:
\begin{eqnarray}
&&\left(\begin{array}{c}  \dot x\\
\dot y\end{array}\right)=\left(\begin{array}{c} ax - by \\ bx+ay\end{array}\right)+M(x,y), \quad  \ \ {\rm if}\ y-\eps<0,\label{neuron1}\\
&&x\to -k\eps,\hskip4.75cm{\rm if}\ y-\eps=0,\label{neuron2}
\end{eqnarray}
where $k>0$, $M(0)=M'(0)=0$, $a<0$, and $b>0$, so that  the origin is a stable focus for subsystem (\ref{neuron1}).

\vskip0.2cm

\noindent In what follows we check the assumptions of Theorem~\ref{thm1}. The impact law (\ref{neuron2}) leads to 
$$
   A'(0)=-k.
$$
The condition (\ref{as2rev}) reduces to
\begin{equation}\label{reduced}
   bu_0(T_0)+a\not=0.
\end{equation}
\noindent To prove the existence of a cycle to the reduced system (\ref{3a})-(\ref{3b}) and to check the condition (\ref{reduced}), we compute $P(A'(0))$ (i.e. $P\left(-k\right)$) for the Poincar\'e map $P$ of  linear system (\ref{3a}) induced by the cross-section $v=-\dfrac{H_\eps(0)}{H_y(0)}=B'(0)=1.$ The linear system (\ref{3a}) corresponding to (\ref{neuron1}) is
\begin{equation}\label{linear1} 
\left(\begin{array}{c}  \dot u\\
\dot v\end{array}\right)=\left(\begin{array}{c} au - bv \\ bu+av\end{array}\right).
\end{equation}
Using that a solution of (\ref{linear1}) is given by
\begin{equation}\label{givenby1}
   u(t)=e^{at}\cos(bt),\quad v(t)=e^{at}\sin(bt),
\end{equation}
we build the following solution of (\ref{linear1})
$$
  u_0(t)=\dfrac{e^{a(t-t_0)}\cos(bt)}{\sin(b t_0)},\quad   v_0(t)=\dfrac{e^{a(t-t_0)}\sin(bt)}{\sin(b t_0)},\quad b t_0={\rm arccot}\left(-k\right),
$$
which verifies the property $(u_0(t_0),v_0(t_0))=\left(-k,1\right)$. It is impossible to find the intersection of solution $(u_0(t),v_0(t))$ with $v=1$ explicitly, so we propose an explicit approach that relies on the observation that 
an intersection of any solution of (\ref{linear1}) with $u=0$ is computable explicitly. 

\vskip0.2cm

\noindent Since ${\rm arccot}\left(-\dfrac{a}{b}\right)\in\left(\dfrac{\pi}{2},\pi\right),$ the first intersection of this solution with $u=0$ occurs at $bt=\dfrac{\pi}{2}+\pi$, which gives
$$
   y_*=v_0\left(\dfrac{1}{b}\cdot\dfrac{3\pi}{2}\right)=-{\exp}\left({a\left(\dfrac{1}{b}\cdot\dfrac{3\pi}{2}-t_0\right)}\right)\dfrac{1}{\sin(b t_0)}.
$$
Now we assume that the intersection of $(u_0(t),v_0(t))$ with $v=1$ occurs at some point $u=r$ and use (\ref{givenby1}) to compute $y_*$ in terms of $r.$ Specifically, using (\ref{givenby1}) we build a solution 
$$
  u^0(t)=\dfrac{e^{a(t-t^0)}\cos(bt)}{\sin(b t^0)},\quad   v^0(t)=\dfrac{e^{a(t-t^0)}\sin(bt)}{\sin(b t^0)},\quad b t^0={\rm arccot}(r),
$$
which verifies $(u^0(t^0),v^0(t^0))=\left(r,1\right).$ Since ${\rm arccot}(r)\in(0,\pi),$ the intersection of $(u^0(t),v^0(t))$ with $u=0,\ v<0,$ must had occurred earlier at time $bt=\dfrac{\pi}{2}-\pi$, which gives
$$
   y^*=v^0\left(\dfrac{1}{b}\cdot\left(-\dfrac{\pi}{2}\right)\right)=-\exp\left(a\left(\dfrac{1}{b}\left(-\dfrac{\pi}{2}\right)-t^0\right)\right)\dfrac{1}{\sin(b t^0)}
$$ 
for the respective point of intersection with $u=0.$ Now equaling $y_*$ and $y^*$, observing that $\dfrac{1}{\sin({\rm arccot}\alpha)}=\sqrt{\alpha^2+1},$ and taking the natural logarithm of both sides of the equality, one gets the following implicit formula for $r:$
\begin{equation}\label{implicit5}
    \dfrac{a}{b}\cdot\dfrac{3\pi}{2}-\dfrac{a}{b}{\rm arccot}\left(-k\right)+\dfrac{1}{2}\ln\left(1+k^2\right)=\psi(r),\   \ \psi(r)=\dfrac{a}{b}\left(-\dfrac{\pi}{2}\right)-\dfrac{a}{b}{\rm arccot}(r)+\dfrac{1}{2}\ln(1+r^2).
\end{equation}
By solving $\psi'(r)=0$ we conclude that $\psi$ is increasing on $r\ge-\dfrac{a}{b}$ and $\psi(r)\to\infty$ as $r\to\infty.$ The equation $\psi(r)=R$ has a solution 
\begin{equation}\label{r0}
r> -\dfrac{a}{b}
\end{equation} for any $R> \psi\left(-\dfrac{a}{b}\right)$. Therefore, $r$ satisfying (\ref{r0}) and (\ref{implicit5}) exists, if
\begin{equation}\label{implicit6}
    \dfrac{a}{b}\cdot\dfrac{3\pi}{2}-\dfrac{a}{b}{\rm arccot}\left(-k\right)+\dfrac{1}{2}\ln\left(1+k^2\right)> \dfrac{a}{b}\left(-\dfrac{\pi}{2}\right)-\dfrac{a}{b}{\rm arccot}\left(-\dfrac{a}{b}\right)+\dfrac{1}{2}\ln\left(1+\dfrac{a^2}{b^2}\right).
\end{equation}
\begin{figure}[h]\center
\includegraphics[scale=0.56]{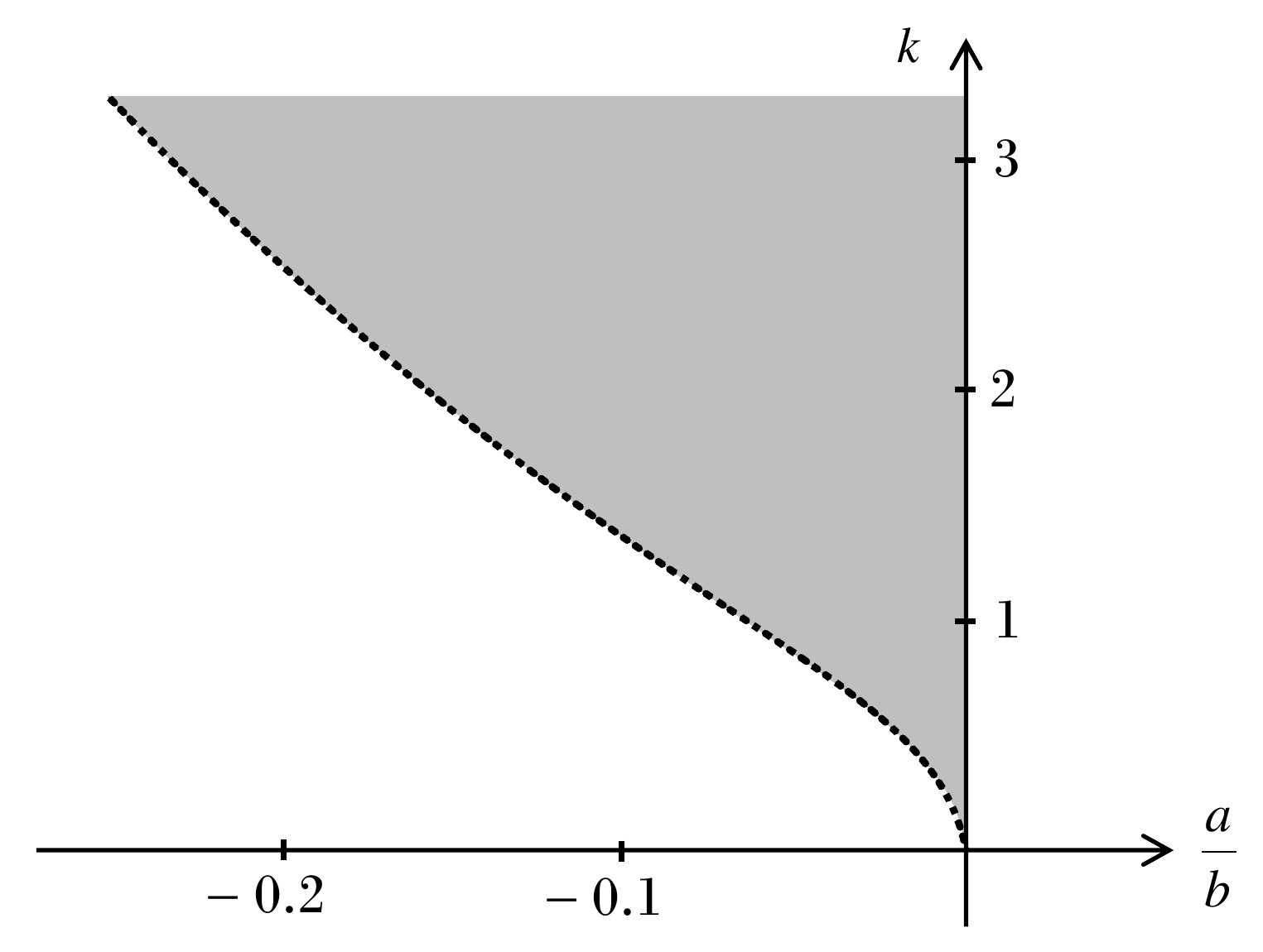}
\vskip-0.2cm
\caption{\footnotesize  The region of parameters $\left(\dfrac{a}{b},k\right)$ that satisfy (\ref{implicit6}), $a<0,$ $b>0$, and $k>0$. The dotted curve is the boundary which is not a part of the region.}  \label{figadd}
\end{figure}
In particular, (\ref{r0}) implies that (\ref{reduced}) holds for the values of $\left(\dfrac{a}{b},k\right)$ satisfying (\ref{implicit6}). 

\vskip0.2cm

\noindent Our findings about the dynamics of (\ref{neuron1})-(\ref{neuron2}) can now be summarized as follows.

\begin{pro}\label{proadd} Assume that $a<0,$ $b>0$, and $k>0$. If $\left(\dfrac{a}{b},k\right)$ satisfies (\ref{implicit6}), then  for all $\eps>0$ sufficiently small, the impacting system (\ref{neuron1})-(\ref{neuron2}) admits a finite-time stable limit cycle $(x_\eps(t),y_\eps(t))$ of one impact per period that shrinks to the origin as $\eps\to 0.$
\end{pro}

\noindent The region of parameters $\left(\dfrac{a}{b},k\right)$ that satisfy the condition of Proposition~\ref{proadd} is plotted in Fig.~\ref{figadd}.

\vskip0.2cm

\noindent Finally, we formulate the following remark that simplifies assumption (\ref{as2rev}) of Theorem~\ref{thm1} in the situations that we are going to consider through the rest of the paper.

\begin{rem}\label{remark} If the impact law (\ref{1b}) satisfies 
\begin{equation}\label{as1rem}
  \left(H_x(A(\eps),B(\eps),\eps),H_y(A(\eps),B(\eps),\eps)\right)\left(\begin{array}{c}f(A(\eps),B(\eps)) \\ g(A(\eps),B(\eps))\end{array}\right)=0,\quad \eps>0,
\end{equation}
then the first assumption of (\ref{as2rev}) reduces to 
 \begin{equation}\label{as2rem}
     u_0(T_0)\not=A'(0).
  \end{equation}
\end{rem}

  %\ \ \mbox{and}\ \  \left(H_x(A(\eps),B(\eps),\eps),H_y(A(\eps),B(\eps),\eps)\right)\left(\begin{array}{c}f(A(\eps),B(\eps)) \\ g(A(\eps),B(\eps))\end{array}\right)=0.

\section{Filippov systems} 

\noindent In this section we consider the following Filippov system equivalent to (\ref{0})
\begin{equation}\label{1}
\left( \begin{array}{c}\dot x\\ \dot y\end{array}\right)=\left( \begin{array}{c}f^i(x,y)\\ g^i(x,y)\end{array}\right),\quad i=\left\{\begin{array}{l}
i=+1,\ {\rm if}\ H(x,y,\eps)>0,\\
i=-1,\ {\rm if}\ H(x,y,\eps)<0,
\end{array}\right.
\end{equation}
where $f^i,$ $g^i$, $i=-1,1$, and $H$ are smooth functions and $\eps>0$ is a parameter.

\begin{pro} \label{pro1} Let  the origin be an equilibrium of the ``$-$''-subsystem of (\ref{1}) and $H(0)=0$. Let  the coordinates be rotated so that 
 $H_x(0)=0$ and $H_y(0)\not=0$. Assume that
 \begin{equation}\label{aspro}
    H_y(0)g^+(0)<0,\quad g_x^-(0)\not=0
 \end{equation}
and
 \begin{equation}\label{excludes}
 \left(\begin{array}{c} f^-_x(0) \\ g^-_x(0)\end{array}\right)\nparallel \left(\begin{array}{c} f^+(0) \\ g^+(0)\end{array}\right).
 \end{equation}
 Then, one can find $r>0$ and $\eps_0>0$ such that for any $0<\eps\le \eps_0$ there exists a unique point  $(A(\eps),B(\eps))\in[-r,r]\times[-r,r]$ which satisfies the property (\ref{as1rem}).  The following properties hold on top of (\ref{as1rem}):
 \begin{itemize}
 \item[1)] The point $(A(\eps),B(\eps))$ satisfies
\begin{equation}\label{A'B'}
   (A'(0),B'(0))=\dfrac{H_\eps(0)}{H_y(0)}\left(\dfrac{g^-_y(0)}{g^-_x(0)},-1\right).
 \end{equation}
 
 \item[2)] The point $(A(\eps),B(\eps))$ splits
$$
    L=\left\{(x,y)\in\mathbb{R}^2: \ x,y\in[-r,r],\ H(x,y,\eps)=0\right\}
 $$ 
 into two parts
$$
     L_{sliding}=\left\{(x,y)\in\mathbb{R}^2: \ x,y\in[-r,r],\ H(x,y,\eps)=0,\ H_{xy}(x,y,\eps)\left(\begin{array}{c}
     f^-(x,y)\\ g^-(x,y)
     \end{array}\right)> 0\right\}
$$
and 
$$
     L_{crossing}=\left\{(x,y)\in\mathbb{R}^2: \ x,y\in[-r,r],\ H(x,y,\eps)=0,\ H_{xy}(x,y,\eps)\left(\begin{array}{c}
     f^-(x,y)\\ g^-(x,y)
     \end{array}\right)< 0\right\}.
$$

 \item[3)] The Filippov equilibrium equation
\begin{eqnarray*}
&& f^-(a,b)-\lambda f^+(a,b)=0,\\
& & g^-(a,b)-\lambda g^+(a,b)=0
\end{eqnarray*}
 possesses a unique equilibrium $(a(\eps),b(\eps),\lambda(\eps))$ on $L$ whose derivative $(a'(0),b'(0),\lambda'(0))$ equals  
 \begin{equation}\label{a'b'}
    \dfrac{H_\eps(0)}{H_y(0)}\left(\dfrac{f^-_y(0)g^+(0)-g^-_y(0)f^+(0)}{f_x^-(0)g^+(0)-g^-_x(0)f^+(0)},-1,-\dfrac{{\rm det}\left\|\begin{array}{cc}  f_x(0) & f_y(0) 
  \\ g_x(0) & g_y(0)\end{array}\right\|}{f_x^-(0)g^+(0)-g^-_x(0)f^+(0)}\right).
 \end{equation}
 
 \item[4)] If 
 \begin{equation}\label{newas}
    \left(f^-_x(0),f^-_y(0)\right)\left(\begin{array}{c} A'(0)\\ B'(0)\end{array}\right) (A'(0)-a'(0))\lambda'(0)<0,
\end{equation} 
% \begin{equation}\label{newas}
%    \dfrac{1}{g_x^-(0)}\det\left\|\begin{array}{cc} f^-_x(0) & f^-_y(0) \\
 %   g^-_x(0) & g^-_y(0) \end{array}\right\|\cdot \left(\dfrac{g^-_y(0)}{g_x^-(0)}-   \dfrac{f^-_y(0)g^+(0)-g^-_y(0)f^+(0)}{-f_x^-(0)g^+(0)+g^-_x(0)f^+(0)}   \right)>0,
% \end{equation}
 then the vector $\left(\begin{array}{c}  f^-(A(\eps),B(\eps)) \\ g^-(A(\eps),B(\eps))  \end{array}\right)$ (tangent to $L$ by definition) points outwards 
   $L_{sliding}.$
 
\noindent \item[5)] If condition (\ref{newas}) holds, then any solution $(x(t),y(t))$ of (\ref{1}) with the initial condition  $(x(0),y(0))$ from the $\left((a(\eps),b(\eps)),(A(\eps),B(\eps))\right)$-segment of  $L_{sliding}$, escapes from  $L_{sliding}$ in finite time through the point $(A(\eps),B(\eps))$. 

\item[6)] The  solution $(x(t),y(t))$ of (\ref{1}) with the initial condition $(x(0),y(0))=(A(\eps),B(\eps))$ leaves $L$ towards 
$$
   L^-=\left\{(x,y)\in\mathbb{R}^2:H(x,y,\eps)<0\right\}
$$
immediately, in the sense that there exists $\Delta t$ such that $t\mapsto (x(t),y(t))$ verifies both the ``$-$''-subsystem of (\ref{1}) and $(x(t),y(t))\in L^-$, for all $t\in(0,\Delta t].$

 \end{itemize}
\end{pro}

\noindent {\bf Proof.}   The existence, uniqueness, and continuous differentiability of $(A(\eps),B(\eps))$ satisfying (\ref{as1rem}) follow by applying the Implicit Function Theorem to the function
$$
   F(A,B,\eps)=\left(\begin{array}{c}H_{xy}(A,B,\eps)\left(\begin{array}{c}f^-(A,B)\\ g^-(A,B)\end{array}\right)\\ H(A,B,\eps)\end{array}\right),
$$
where we use that   $F(0)=0$ and $\det\|F_{AB}(0)\|\not=0$ by  the second of the assumptions of (\ref{aspro}).

\vskip0.2cm

\noindent {\bf Part 1.} Formula (\ref{A'B'}) follows by computing the derivative of $F(A(\eps),B(\eps),\eps)=0$ at 
$\eps=0.$

\vskip0.2cm

\noindent {\bf Part 2.} Follows from the uniqueness of $(A(\eps),B(\eps))$.

\vskip0.2cm

\noindent {\bf Part 3.} The region $L_{sliding}$ is the region of sliding by the first of the assumptions of (\ref{aspro}). We define $(a(\eps),b(\eps))$ as the unique equilibrium of the sliding vector field of Filippov system (\ref{1}). To prove the existence of such a unique equilibrium we apply the Implicit Function Theorem to the function
$$
   G(a,b,\lambda,\eps)=\left(\begin{array}{c}
   f^-(a,b)-\lambda f^+(a,b)\\
   g^-(a,b)-\lambda g^+(a,b) \\ H(a,b,\eps)\end{array}\right).
$$
The determinant 
$$
{\rm det}|G_{ab\lambda}(0)|={\rm det}\left|\begin{array}{ccc}
f^-_x(0) & f_y^-(0) & -f^+(0) \\
g^-_x(0) & g^-_y(0) & -g^+(0)\\
0 & H_y(0) & 0 
\end{array}\right|=-H_y(0)(-f_x^-(0)g^+(0)+g_x^-(0)f^+(0))
$$ doesn't vanish by (\ref{excludes}) and the formula for the derivative of the implicit function
\begin{equation}\label{formula}
   (a'(0),b'(0),\lambda'(0))^T=-G_{ab\lambda}(0)^{-1}G_\eps(0)
\end{equation}
 yields (\ref{a'b'}). 

\vskip0.2cm

\noindent {\bf Part 4.}  Conditions (\ref{A'B'}) and (\ref{a'b'}) imply that $A(\eps)a(\eps)\not=0$ for all $\eps>0$ sufficiently small. Case I: $\lambda'(0)<0,$ which combined with (\ref{newas}) gives
\begin{equation}\label{newas1}
  \left(f^-_x(0),f^-_y(0)\right)\left(\begin{array}{c} A'(0)\\ B'(0)\end{array}\right) (A'(0)-a'(0))>0.
\end{equation}
Furthermore, $\lambda'(0)<0$ implies that $(a(\eps),b(\eps))\in L_{sliding}$ for all $\eps>0$ sufficiently small. 
Sub-case~1: $A'(0)<a'(0)$ (i.e. $(A(\eps),B(\eps))$ is the left endpoint of $L_{sliding}$). In this case (\ref{newas1}) yields $f^-(A(\eps),B(\eps))<0,$ i.e. the vector $\left(\begin{array}{c}  f^-(A(\eps),B(\eps)) \\ g^-(A(\eps),B(\eps))  \end{array}\right)$ points to the left. Sub-case 2: By analogy, when $A'(0)>a'(0),$ the assumption (\ref{newas1}) implies $f^-(A(\eps),B(\eps))<0.$

\vskip0.2cm

\noindent Case II: $\lambda'(0)>0.$ Can be considered by analogy taking into account that $\lambda'(0)>0$ implies that $(a(\eps),b(\eps))\in L_{crossing}$ for all $\eps>0$ sufficiently small.

\vskip0.2cm

\noindent {\bf Part 5.} The dynamics of $(x(t),y(t))$ is described by one-dimensional smooth equation of sliding motion (Filippov~\cite[\S19]{filippov}) as long as $(x(t),y(t))\in L_{sliding}.$ 
Part 4) implies that the vector field of the equation of sliding motion on $L_{sliding}$ points towards the endpoint $(A(\eps),B(\eps))$ at all the points of $L_{sliding}$ close to $(A(\eps),B(\eps))$. 
Therefore, if we assume, by contradiction, that the solution $(x(t),y(t))$ doesn't reach $(A(\eps),B(\eps))$ in finite-time, then the sliding vector field must possess an equilibrium on the $\left((a(\eps),b(\eps)),(A(\eps),B(\eps))\right)$-segment of  $L_{sliding}$, which contradicts the uniqueness of equilibrium $(a(\eps),b(\eps)).$

\vskip0.2cm

\noindent {\bf Part 6.} This is a standard property, see e.g. Filippov~\cite[\S19]{filippov}.

\vskip0.2cm

\noindent The proof of the proposition is complete. \qed

\vskip0.3cm

\noindent Combining Theorem~\ref{thm1} (where we view the ``$-$''-subsystem of (\ref{1}) as system (\ref{1a})), Remark~\ref{remark}, and Proposition~\ref{pro1}, we arrive to the following result about limit cycles of Filippov system (\ref{1}).

\begin{theorem}\label{thm2}  Let  the origin be an equilibrium of the ``$-$''-subsystem of (\ref{1}) and $H(0)=0$. Let  the coordinates be rotated so that 
 $H_x(0)=0$ and $H_y(0)\not=0$. Let the assumptions
 (\ref{aspro}), 
  (\ref{excludes}), and (\ref{newas}) of Proposition~\ref{pro1} hold with $(A'(0),B'(0))$ and $(a'(0),b'(0),\lambda'(0))$ given by (\ref{A'B'}) and (\ref{a'b'}) respectively. 
Assume that the reduced system (\ref{3a})-(\ref{3b}) with $(f,g)$ replaced by $(f^-,g^-)$ admits a cycle $(u_0(t),v_0(t))$  with the initial condition $(u_0(0),v_0(0))=(A'(0),B'(0))$
of exactly one impact per period. Let $T_0$ be the period of the cycle. If 
\begin{equation}\label{segment}
\begin{array}{l}
  u_0(T_0)\in \left(\min\{a'(0),A'(0)\},\max\{a'(0),A'(0)\}\right)\quad\mbox{in the case when }\lambda'(0)<0,\\
  u_0(T_0)\not=A'(0)\quad\mbox{in the case when }\lambda'(0)>0,
  \end{array}
\end{equation}
then for all $\eps>0$ sufficiently small, the Filippov system (\ref{1}) admits a finite-time stable stick-slip limit cycle $(x_\eps(t),y_\eps(t))\to 0$  as $\eps\to 0.$ 
\end{theorem}

\noindent {\bf Proof.} Let $(x_\eps(t),y_\eps(t))$ be the solution of (\ref{1}) with the initial condition $(x_\eps(0),y_\eps(0))=(A(\eps),B(\eps))$ as defined in Theorem~\ref{thm1}. To prove the theorem it is sufficient to observe that condition (\ref{segment}) implies that $(x_\eps(T_\eps),y_\eps(T_\eps))$ belongs to t the  $\left((a(\eps),b(\eps)),(A(\eps),B(\eps))\right)$-segment of  $L_{sliding}$ as defined in Proposition~\ref{pro1}, so that the map (\ref{1b}) is well defined on $L$ in the neighborhood of $(x_\eps(T_\eps),y_\eps(T_\eps)).$ \qed

\vskip0.2cm

\noindent As an example, we consider the following Filippov system
\begin{eqnarray}
&& \left(\begin{array}{c}  \dot x\\
\dot y\end{array}\right)=\left(\begin{array}{c} m \\ -1\end{array}\right)+K(x,y), \hskip1.45cm {\rm if}\ y-\eps>0,\label{nform1}\\
&& \left(\begin{array}{c}  \dot x\\
\dot y\end{array}\right)=\left(\begin{array}{c} ax - by \\ bx+ay\end{array}\right)+M(x,y), \quad  \ \ {\rm if}\ y-\eps<0,\label{nform2}
\end{eqnarray}
where $a,b>0,$ $m\in\mathbb{R}$, and $K,M$ are any $C^2$ function such that $K(0)=M'(0)=0.$ 

\vskip0.2cm

\noindent In what follows we check the assumptions of Theorem~\ref{thm2}. Assumptions (\ref{aspro}) and (\ref{excludes}) hold, if $a\not=0$ and $\dfrac{a}{b}\not=-m$ respectively.  Formulas (\ref{A'B'}) and (\ref{a'b'}) lead to the following expressions for the derivatives $A'(0),$ $B'(0),$ $a'(0),$ $b'(0),$ and $\lambda'(0):$
\begin{eqnarray}
  (A'(0),B'(0)) &=& \left(-\dfrac{a}{b},1\right),\nonumber\\
  (a'(0),b'(0),\lambda'(0))&=& \left(\dfrac{b-am}{a+bm},1,-\dfrac{a^2+b^2}{a+bm}\right),\label{lambda}
\end{eqnarray}
which gives
$$
   \dfrac{a^2+b^2}{b}\cdot\dfrac{-a^2-b^2}{b(a+bm)}\cdot\dfrac{a^2+b^2}{a+bm}
$$
for the left-hand-side of (\ref{newas}). Therefore, assumption (\ref{newas}) always holds.

\vskip0.2cm

\noindent To prove the existence of a cycle to the reduced system (\ref{3a})-(\ref{3b}) and to check the condition (\ref{segment}), we have to compute $r=P(A'(0))=P\left(-\dfrac{a}{b}\right)$. But the same quantity $r=P\left(-k\right)$ has been already computed in the example of Section~2. Therefore, to obtain the formula for $r$ we simply need to replace $k$ by $\dfrac{a}{b}$ in formula (\ref{implicit5}) of Section~2 getting
\begin{equation}\label{implicit}
    \dfrac{a}{b}\cdot\dfrac{3\pi}{2}-\dfrac{a}{b}{\rm arccot}\left(-\dfrac{a}{b}\right)+\dfrac{1}{2}\ln\left(1+\dfrac{a^2}{b^2}\right)=\dfrac{a}{b}\left(-\dfrac{\pi}{2}\right)-\dfrac{a}{b}{\rm arccot}(r)+\dfrac{1}{2}\ln(1+r^2).
\end{equation}
\begin{figure}[h]\center
\includegraphics[scale=0.5]{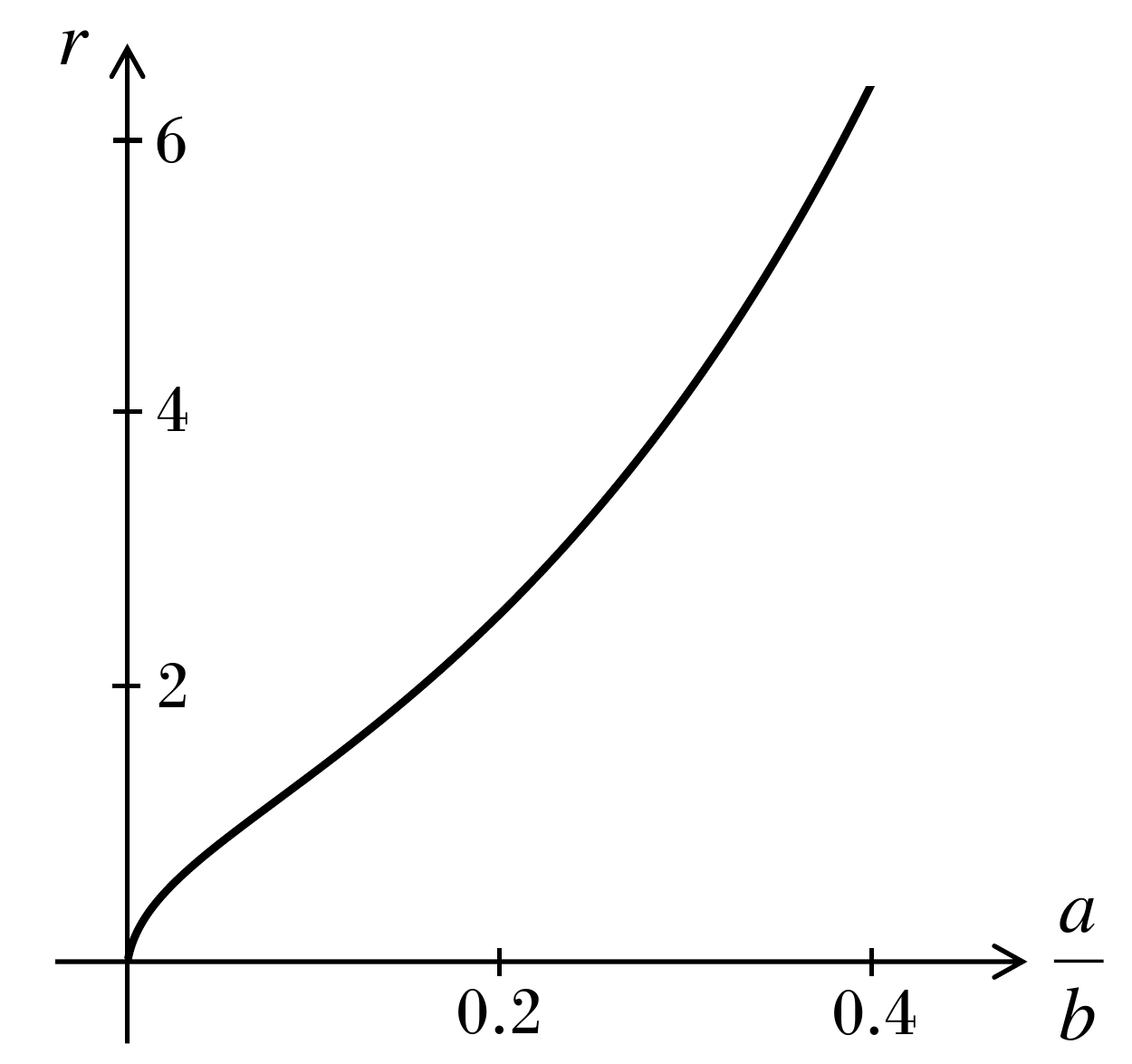}\qquad\qquad  \includegraphics[scale=0.56]{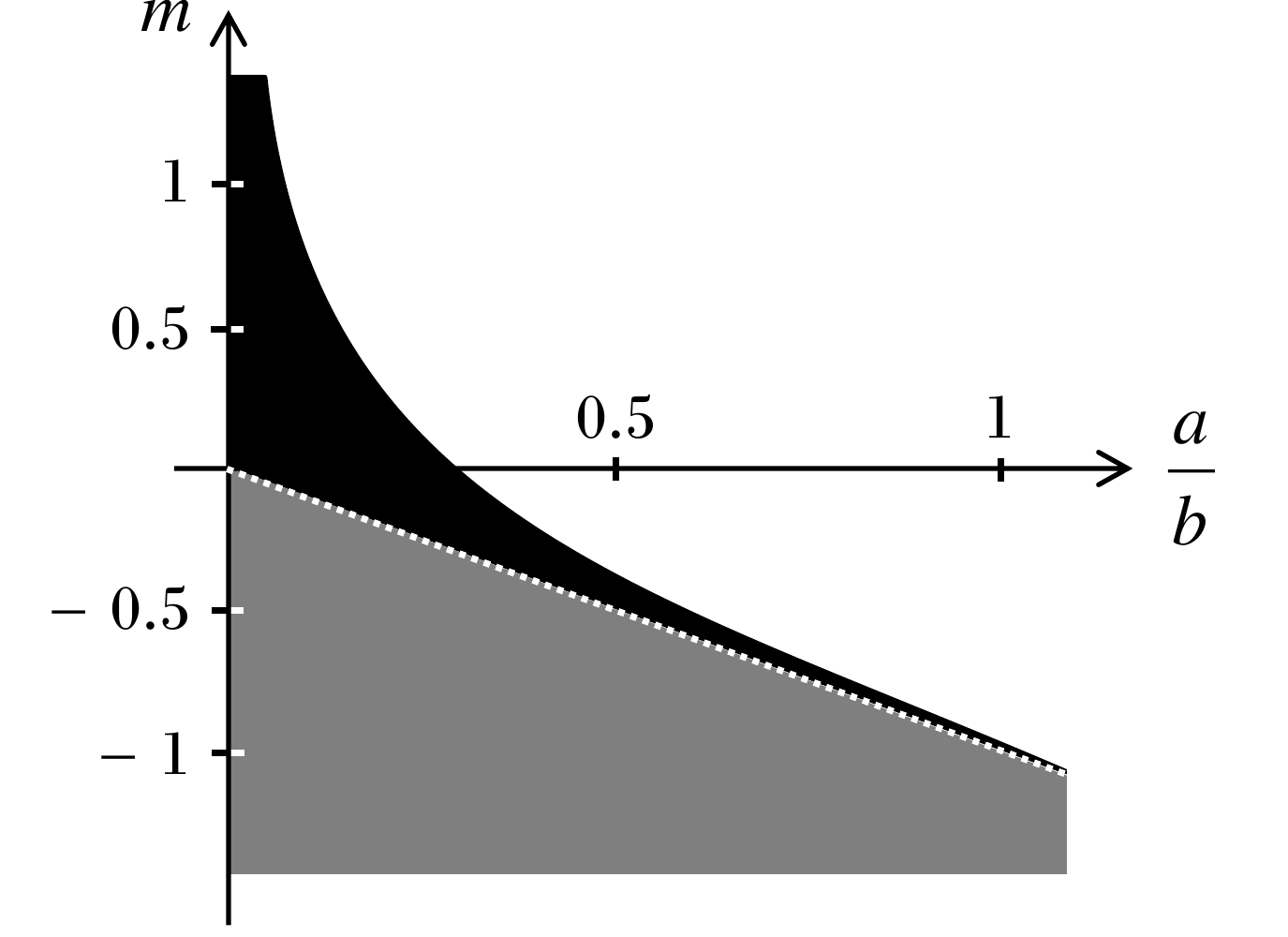}
\vskip-0.2cm
\caption{\footnotesize Left: The solution of (\ref{implicit}). Right:  The region of parameters $\left(\dfrac{a}{b},m\right)$ that satisfy (\ref{m1}) (gray), the region of parameters $\left(\dfrac{a}{b},m\right)$ that satisfy  (\ref{m2})-(\ref{implicit1}) (black), and the line $m=-\dfrac{a}{b}$ (dotted white).} \label{fig1}
\end{figure}
The graph of the implicit equation (\ref{implicit}) is given in Fig.~\ref{fig1}~left, from which we conclude that the solution $(u_0(t),v_0(t))$ returns back to the cross-section $v=1$ at the value $r\left(\dfrac{a}{b}\right)=(u_0(T_0),v_0(T_0))$ which increases monotonically with  $\dfrac{a}{b}.$ To summarize, the requirement of Theorem~\ref{thm2} about the existence of a cycle to the reduced system (\ref{3a})-(\ref{3b}) holds. Our goal now is to establish (\ref{segment}).

\vskip0.2cm

\noindent Based on (\ref{lambda}), the property $\lambda'(0)>0$ is equivalent to 
\begin{equation}\label{m1}
m<-\dfrac{a}{b}.
\end{equation}
Therefore, if (\ref{m1}) is satisfied, then the assumption (\ref{segment}) of Theorem~\ref{thm2} holds. Let us consider  $\lambda'(0)<0,$ i.e.
\begin{equation}\label{m2}
m>-\dfrac{a}{b}.
\end{equation}
In this case assumption (\ref{segment}) takes the form
$$   r=u_0(T_0)<\dfrac{1-\dfrac{a}{b}m}{\dfrac{a}{b}+m}.
$$
Since $r\mapsto -\dfrac{a}{b}{\rm arccot}(r)+\dfrac{1}{2}\ln(1+r^2)$ is a monotonically increasing function, we can combine the later inequality with (\ref{implicit}) to obtain
\begin{equation}\label{implicit1}
\begin{array}{l}
    \dfrac{a}{b}\cdot\dfrac{3\pi}{2}-\dfrac{a}{b}{\rm arccot}\left(-\dfrac{a}{b}\right)+\dfrac{1}{2}\ln\left(1+\dfrac{a^2}{b^2}\right)<\\\qquad <\dfrac{a}{b}\left(-\dfrac{\pi}{2}\right)-\dfrac{a}{b}{\rm arccot}\left(\dfrac{1-\dfrac{a}{b}m}{\dfrac{a}{b}+m}\right)+\dfrac{1}{2}\ln\left(1+\left(\dfrac{1-\dfrac{a}{b}m}{\dfrac{a}{b}+m}\right)^2\right).
    \end{array}
\end{equation}
We arrive to the following corollary of Theorem~\ref{thm2}.
\begin{pro} If $\left(\dfrac{a}{b},m\right)$ satisfies either (\ref{m1}) or   (\ref{m2})-(\ref{implicit1}), then  for all $\eps>0$ sufficiently small, the Filippov system (\ref{nform1})-(\ref{nform2}) admits a finite-time stable stick-slip limit cycle $(x_\eps(t),y_\eps(t))$ that shrinks to the origin as $\eps\to 0.$
\end{pro}

\noindent The region of parameters $\left(\dfrac{a}{b},m\right)$ that satisfy (\ref{m1}) and the region of parameters $\left(\dfrac{a}{b},m\right)$ that satisfy  (\ref{m2})-(\ref{implicit1}) are drawn at Fig.~\ref{fig1}~right.

\section{Sweeping processes}

Consider a perturbed sweeping process
\begin{equation}\label{sp}
   \left(\begin{array}{c}\dot x\\ \dot y\end{array}\right)\in -N_{C(\eps)}(x,y)+\left(\begin{array}{c}  f(x,y) \\ g(x,y)\end{array}\right),
\end{equation}
where 
$$
   C(\eps)=\left\{(x,y)\in\mathbb{R}^2:H(x,y,\eps)\le 0\right\}
$$
is a nonempty convex closed time-independent set, for all $\eps\ge 0.$ We will assume that $f$ and $g$ are $C^1$ globally Lipschitz functions, so that for any initial condition $(x_0,y_0)\in C(\eps)$, the sweeping process (\ref{sp}) admits a unique forward solution $(x(t),y(t))\in C(\eps)$ with the initial condition $(x(0),y(0))=(x_0,y_0)$, that satisfies the differential inclusion (\ref{sp}) for a.a. $t\ge 0$ (Edmond-Thibault \cite[Theorem~1]{edm}). According to the definition of the solution $(x(t),y(t))$, 
\begin{equation}\label{np}
  \left(\begin{array}{c}\dot x\\ \dot y\end{array}\right)=\left(\begin{array}{c}  f(x,y) \\ g(x,y)\end{array}\right),\quad \mbox{when}\ (x(t),y(t))\in {\rm int}C(t).
\end{equation}

\begin{pro} \label{pro2} Let  the origin be an equilibrium of the subsystem  (\ref{np}) and $H(0)=0$. Let  the coordinates be rotated so that 
 $H_x(0)=0$ and $H_y(0)\not=0$. Assume that
 \begin{equation}\label{excludesp}
 f_x(0)\not=0,\quad g_x(0)\not=0.
 \end{equation}
 Then, there exist $r>0$ and $\eps_0>0$ such that for any $0<\eps\le \eps_0$ there exists a unique point $(A(\eps),B(\eps))\in[-r,r]\times[-r,r]$ which satisfies the property (\ref{as1}).  The statements 1)-2) and 4)-6) of Proposition~\ref{pro1} hold with $f^-$, $g^-$, and (\ref{1}) replaced by $f$, $g$, and (\ref{sp}) respectively. Furthermore, the following analogue of statement 3) of Proposition~\ref{pro1} takes place
 \begin{itemize}
  \item[3)] a) Any solution $(x(t),y(t))$ of sweeping process (\ref{1}) with the initial condition  $(x(0),y(0))\in L_{sliding}$ can  escape from  $L_{sliding}$ through the endpoints of $L_{sliding}$ only (i.e. through the two points of $\overline{L_{sliding}}\backslash L_{sliding}$). b) While in $L_{sliding}$, the solution $(x(t),y(t))$ is governed by the following equation of sliding motion  
\begin{equation}\label{ssp}
  \left(\begin{array}{c} \dot x(t) \\ \dot y(t)\end{array}\right)=\dfrac{ 
   \left<\left(\begin{array}{c} f(x(t),y(t)) \\ g(x(t),y(t)) \end{array}\right),\left(\begin{array}{c}-H_y(x(t),y(t),\eps) \\  H_x(x(t),y(t),\eps)\end{array}\right)\right>}{\|H_{xy}(x(t),y(t),\eps)\|}\left(\begin{array}{c}-H_y(x(t),y(t),\eps) \\  H_x(x(t),y(t),\eps)\end{array}\right).  
\end{equation}
c) The equation
\begin{equation}\label{construction}
\begin{array}{l}
 f^-(a,b)+\lambda H_x(a,b,\eps)=0,\\
 g^-(a,b)+\lambda H_y(a,b,\eps)=0
\end{array}
\end{equation}
for the equilibrium of (\ref{ssp})
  possesses a unique solution $(a(\eps),b(\eps),\lambda(\eps))$ on $L$ with
 \begin{equation}\label{a'b'sp}
    (a'(0),b'(0),\lambda'(0))=\dfrac{H_\eps(0)}{H_y(0)}\left(\dfrac{f_y(0)}{f_x(0)},-1,\dfrac{1}{H_y(0)f_x(0)}{\rm det}\left|\begin{array}{cc}
    f_x(0) & f_y(0) \\ g_x(0) & g_y(0)
    \end{array}\right|\right).
 \end{equation}
\end{itemize}
\end{pro}

\noindent{\bf Proof.} {\bf Part 1 and Part 2.} Same proof as in Proposition~\ref{pro1}, where the second of the assumptions of (\ref{excludesp}) is used.

\vskip0.2cm

\noindent {\bf Part 3a.}  Fix $\eps>0$. Let $t_{escape}\ge 0$ be the time when $(x(t),y(t))$ escapes from $L_{sliding}$, specifically
$$
    t_{escape}=\max\{t_0\ge 0:x(t)\in[-r,r],\ y(t)\in[-r,r],\ H(x(t),y(t),\eps)=0,\ t\in[0,t_0]\}.
$$
Assuming that neither $|x(t_{escape})|=r,$ nor $|y(t_{escape})|=r,$ we now show that
\begin{equation}\label{nowshow}
H_{xy}(x(t_{escape}),y(t_{escape}),\eps)\left(\begin{array}{c}
     f(x(t_{escape}),y(t_{escape}))\\ g(x(t_{escape}),y(t_{escape}))
     \end{array}\right)\le 0,
\end{equation}
which coincides with the Statement 3a.

\vskip0.2cm

\noindent By the definition of $t_{escape},$ for any $\delta>0$ there exist $t_\delta\in[t_{escape},t_{escape}+\delta]$ such that $H(x(t_\delta),y(t_\delta),\eps)<0$ and $t_\delta^*\in[t_{escape},t_\delta)$ such that 
$$
   H(x(t_\delta^*),y(t_\delta^*),\eps)=0,\quad H(x(t),y(t),\eps)\not=0,\quad t\in(t_\delta^*,t_\delta].
$$
Since,  the solution $(x(t),y(t))$ satisfies (\ref{np}) on $(t_\delta^*,t_\delta]$, one can apply the Mean-Value Theorem to get
$$
   \left(\begin{array}{c}
     x(t_\delta) \\
     y(t_\delta)
     \end{array}\right)=\left(\begin{array}{c} x(t_\delta^*) \\ y(t_\delta^*) \end{array}\right)
     + \left(\begin{array}{c}  f(x(t_\delta^{**}),y(t_\delta^{**})) \\
     g(x(t_\delta^{**}),y(t_\delta^{**})) \end{array}\right)(t_\delta-t_\delta^*), 
$$ 
or
$$
   H(x(t_\delta^*) + f(x(t_\delta^{**}),y(t_\delta^{**})),y(t_\delta^*) + g(x(t_\delta^{**}),y(t_\delta^{**})),\eps)<0,
$$
which yields (\ref{nowshow}) as $\delta\to 0.$

\vskip0.2cm

\noindent {\bf Part 3b.} Consider some $t_0>0$ such that $(x(t),y(t))\in L_{sliding}$ for all $t\in[0,t_0].$ From the definition of $L_{sliding}$ we conclude that 
$$
\left<\left(\begin{array}{c} \dot x(t) \\ \dot y(t) \end{array}\right),H_{xy}(x(t),y(t),\eps)\right>=0,\quad \mbox{for a.a.}\ t\in[0,t_0],
$$
where we use that the derivatives of solutions of (\ref{sp}) are defined for a.a. $t$ only, and so
$$
\dfrac{1}{\|H_{xy}(x(t),y(t),\eps)\|}\left<\left(\begin{array}{c} \dot x(t) \\ \dot y(t) \end{array}\right),\left(\begin{array}{c}-H_y(x(t),y(t),\eps) \\  H_x(x(t),y(t),\eps)\end{array}\right)\right>\left(\begin{array}{c}-H_y(x(t),y(t),\eps) \\  H_x(x(t),y(t),\eps)\end{array}\right)=\left(\begin{array}{c} \dot x(t) \\ \dot y(t) \end{array}\right),
$$
for a.a. $t\in[0,t_0].$ Equation (\ref{ssp}) now comes by projecting (\ref{sp}) on the vector  
$\left(\begin{array}{c}-H_y(x(t),y(t),\eps) \\  H_x(x(t),y(t),\eps)\end{array}\right)$ and by extending (\ref{ssp}) from a.a. $t\in[0,t_0]$ to all $t\in[0,t_0]$ using smoothness of (\ref{ssp}).

\vskip0.2cm

\noindent {\bf Part 3c.} To prove the existence and uniqueness of $(a(\eps),b(\eps))$, we apply the Implicit Function Theorem to the function
$$
   G(a,b,\lambda,\eps)=\left(\begin{array}{c}
   f(a,b)+\lambda H_x(a,b,\eps)\\
   g(a,b)+\lambda H_y(a,b,\eps) \\ H(a,b,\eps)\end{array}\right).
$$
The determinant 
$$
{\rm det}|G_{ab\lambda}(0)|={\rm det}\left|\begin{array}{ccc}
f^-_x(0) & f_y^-(0) & 0 \\
g^-_x(0) & g^-_y(0) & H_y(0)\\
0 & H_y(0) & 0 
\end{array}\right|=-H_y(0)^2f_x^-(0)
$$ doesn't vanish by the first assumption of (\ref{excludesp}) and the formula (\ref{formula}) for the derivative of the implicit function
 yields (\ref{a'b'sp}).

\vskip0.2cm

\noindent {\bf Part 4 and Part 5.} Same proof as in Proposition~\ref{pro1}. In particular, the construction (\ref{construction}) implies that $(a(\eps),b(\eps))\in L_{sliding}$ for all $\eps>0$ sufficiently small, if $\lambda'(0)<0$, and $(a(\eps),b(\eps))\in L_{crossing}$ for all $\eps>0$ sufficiently small, if $\lambda'(0)>0$.

 \vskip0.2cm
 
  \noindent {\bf Part 6.}  Let $(x(t),y(t))$ be the solution of (\ref{np}) with the initial condition $(x(0),y(0))=(A(\eps),B(\eps)).$ By the definition of $(A(\eps),B(\eps))$, there exists $\Delta t>0$ such that $H(x(t),y(t),\eps)<0$ for all $t\in(0,\Delta t].$ Therefore, $(x(t),y(t))$ is the solution of (\ref{sp}) on $(0,\Delta t].$ Therefore, $(x(t),y(t))$ is the solution of (\ref{sp}) on $[0,\Delta t],$ because the definition of the solution (\ref{sp}) requires the validity of (\ref{sp}) for $(x(t),y(t))$ in a.a. time instances $t$ only.
  
  \vskip0.2cm
  
  \noindent The proof of the proposition is complete.
  \qed

\vskip0.3cm

\noindent Combining Theorem~\ref{thm1} (where we view the ``$-$''-subsystem of (\ref{1}) as system (\ref{1a})) and Proposition~\ref{pro1}, we arrive to the following result about limit cycles of Filippov system (\ref{1}).

\begin{theorem}\label{thm2}  Let  the origin be an equilibrium of the subsystem of (\ref{np}) and $H(0)=0$. Assume that  the coordinates are rotated so that 
 $H_x(0)=0$ and $H_y(0)\not=0$. Let the assumption
   (\ref{newas}) of Proposition~\ref{pro1} hold with $(A'(0),B'(0))$ and $(a'(0),b'(0),\lambda'(0))$ given by (\ref{A'B'}) and (\ref{a'b'sp}) respectively. Let the assumption (\ref{excludesp}) of Proposition~\ref{pro2} holds.
Finally, assume that the reduced system (\ref{3a})-(\ref{3b}) admits a cycle $(u_0(t),v_0(t))$  with the initial condition $(u_0(0),v_0(0))=(A'(0),B'(0))$
of exactly one impact per period. Let $T_0$ be the period of the cycle. If (\ref{segment}) holds then for all $\eps>0$ sufficiently small, the sweeping process (\ref{sp}) admits a finite-time stable stick-slip limit cycle $(x_\eps(t),y_\eps(t))\to 0$  as $\eps\to 0.$ 
\end{theorem}

\begin{figure}[h]\center
\includegraphics[scale=0.5]{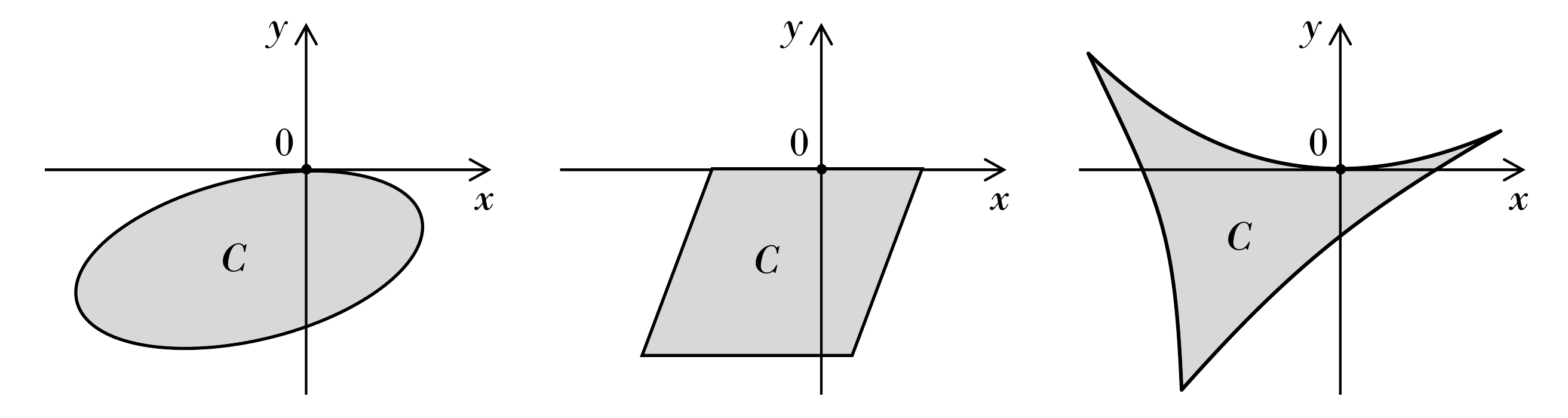}\qquad\qquad  \vskip-0.2cm
\caption{\footnotesize Two convex sets and an $r$-prox regular set that can be used in sweeping process (\ref{spex}).} \label{fig2}
\end{figure}

\noindent To illustrate the theorem we will build upon computations from the example of Section~3 and consider the following sweeping process
\begin{equation}\label{spex}
   \left(\begin{array}{c}\dot x\\ \dot y\end{array}\right)\in -N_{C-\eps\left(\begin{array}{c}0\\ 1\end{array}\right)}(x,y)+\left(\begin{array}{c} ax - by \\ bx+ay\end{array}\right)+M(x,y),
\end{equation}
where $C$ is any compact convex set or an $r$-prox regular whose boundary $\partial C$ contains the origin and 
$$
   \partial C=\left\{(x,y)\in\mathbb{R}^2:H(x,y)=0\right\}\quad\mbox{in  the neighborhood of the origin,}
$$
where $H$ is $C^1$-function such that $H_{xy}(0,0)=\left(\begin{array}{c} 0 \\ 1\end{array}\right),$ see Fig.~\ref{fig2}.

\noindent In order to adopt computations of the Example of Section~3 we only have to replace $(a'(0),b'(0),\lambda'(0))$ of (\ref{a'b'}) by $(a'(0),b'(0),\lambda'(0))$ of (\ref{a'b'sp}) when computing estimates 
(\ref{m1})-(\ref{implicit1}). From  (\ref{a'b'sp})  we have
$$
   \left(\dfrac{b}{a},1,-\dfrac{a^2+b^2}{a}\right),
$$
which just equals (\ref{lambda}) with $m=0.$ The next proposition, therefore, comes by plugging $m=0$ into (\ref{implicit1}).
\begin{pro}
 If $\dfrac{a}{b}$ satisfies 
 $$
\dfrac{a}{b}\left(4{\rm arctan}\dfrac{a}{b}-3\pi\right)>2\ln\dfrac{a}{b}\qquad \mbox{(which gives approximately }\dfrac{a}{b}<0.29),
 $$
 then  for all $\eps>0$ sufficiently small, the sweeping process  (\ref{spex}) admits a finite-time stable stick-slip limit cycle $(x_\eps(t),y_\eps(t))$ that shrinks to the origin as $\eps\to 0.$
\end{pro}

\section{Conclusions} The results of this paper complement the available literature in various ways. First of all, our theorem on bifurcation of limit cycles from a boundary equilibrium of an impacting system turned out to be applicable in the case of a stable focus, thus giving a proof for the occurrence of spiking oscillations in a simple resonate-and-fire model. 

\vskip0.2cm

\noindent Even though studies on bifurcation of limit cycle from a focus boundary equilibrium in Filippov systems are extensively available, our approach establishes the occurrence of limit cycles in the initial Filippov system, rather than in its reduced normal form.

\vskip0.2cm

\noindent Perhaps most importantly, this paper offers the first ever result on bifurcation of limit cycles in sweeping processes, in which analysis we derived an equation of sliding along the boundary of an unilateral constraint and observed that the action of the unilateral constraint is equivalent to an action of an orthogonal vector field pointing towards the unilateral constraint from the outside.

\end{document}